\documentclass{llncs}

\usepackage{makeidx}  
\usepackage{amsfonts, amsmath, amsfonts, amssymb}
\usepackage{graphicx} %
\usepackage{tikz}

\begin{document}
\frontmatter          
\pagestyle{headings}

\title{On Lyndon-Word Representable Graphs}
\titlerunning{Lyndon-Word Representable}  
%
\author{Hossein Teimoori Faal}
\authorrunning{H. Teimoori}   
%
\tocauthor{Hossein Teimoori Faal}
\institute{Department of Mathematics and Computer Science, 
Allameh Tabataba'i University, Tehran, Iran,\\
\email{hossein.teimoori@gmail.com}}

\maketitle              

\begin{abstract}

In this short note, we first associate a new simple undirected graph with a 
given word over an ordered alphabet of $n$-letters. We will 
call it the Lyndon graph of that word. Then, we introduce the 
concept of the Lyndon-word representable graph as a graph
isomorphic to a Lyndon graph of some word. Then, we
introduce the generalized Stirling cycle number $S(N;n,k)$
as the number words of length $N$ 
with $k$ distinct Lydon words in their Lyndon
factorization
over an ordered alphabet of 
$n$-letters . 
Finally, we conclude the paper with several interesting 
open questions and conjectures for interested audiences. 

\end{abstract}
\section{Introduction} 
\label{saucan-sec:intro}

A functional digraph of a permutation is a simple 
way to associate a graph with an algebraic structure. 
The authors in \cite{BM2016} have associated 
a simple, finite and undirected graph $G(w)$ by 
a word $w$ (a permutation of a multiset of letters). 
For a given graph $G$, if there exits a word 
$w$ such that its corresponding graph is 
isomorphic to $G(w)$, then the graph $G$ is called a \emph{Parikh-word representable graph}. 
\\
Here, we associate another simple graph with a given
word in a way that the corresponding graph of 
any Lyndon word is always a \emph{connected}
graph. Therefore, we will call it the 
Lyndon graph of a word. 
\\
We will also call a graph
$G$ 
a \emph{Lyndon-word representable graph} if there exists 
a word $w$ such that its associated Lyndon graph $G(w)$ is 
isomorphic to $G$. 
\\
The main purpose of this talk is to give a 
\emph{characterization} of Lydon-word representable 
graphs and their connection with combinatorial 
numbers like the generalized Stirling cycle 
numbers. 

%
\section{Basic Definitions and Notations} 
\label{Definitions}

In analogy with the idea of \emph{word representable graphs}
introduced in \cite{BM2016}, we have the following definiton.

\begin{definition}  \label{firstdef:Lyndon}

Let 
$
\mathcal{A} = 
\{
a_{1} \prec a_{2} \prec \cdots \prec a_{n}
\}
$
be an ordered alphabet of $n$-letters and 
$
w = w_{1}w_{2} \cdots w_{N}
$
be a word of length $N(\geq n)$ over 
$
\mathcal{A}
$
.
We define a \emph{Lyndon graph} 
$
G(w)
$
with the vertex set 
$
[N]
$
and there exists and edge between 
two vertices $i$ and $j$ 
($1 \leq i,j \leq N$) if $w_{i}w_{j}$ is a
scattered subword of $w$ in the form 
$
a_{r}a_{s}
$
for some $r$ and $s$ ($1 \leq r < s \leq n$). 
A graph $G$ is called a \emph{Lyndon-word representable graph}
if there exists a word $w \in \mathcal{A}^{\star}$ such that 
$G$ is isomorphic to $G(w)$.

\end{definition}

\begin{remark}
	It is worth to note that for the alphabet of $2$-letters, 
	the \emph{Parikh-word representable graphs} and the
	{Lyndon-word representable graphs} are the same.
	Thus, one can consider the class of Lyndon-word representable graphs 
	as a generalization of the class of Parikh-word representable graphs. 
	 
\end{remark}

Form now on, we will denote the number of occurences of a letter
$a$ in a word $w$ by $\vert w \vert_{a}$. 
\begin{remark}
	It is note worthy that the word corresponds to a connected 
	Lyndon-word representable graph must start with a letter $a$ (or $b$) and end with a letter $c$. 
	In the case starting with a letter $b$, the $\vert w \vert_{b}$
	must be greater than $1$ and it ends with a letter $c$. 
	It is not hard to see that for a Lyndon-word representable graph
	$G=G(w)$, in which the word $w$ starts with $a$ and ends with $c$
	there are two adjacent vertices in $G$ with labels $1$
	and $N$, respectively such that their degree sum is equal to 
	$N$. 
\end{remark}

\begin{example}
	Let
	$
	\mathcal{A} = \{a \prec b
	\prec c
	\}
	$
	be an ordered ternary alphabet and the word
	$
	w = abccab
	$
	, 
	the graph $G(w)$ has depicted in Figure \ref{fig:results}.
	
	\begin{figure}[h]
		\begin{center}
			\begin{tikzpicture}
				\draw (0,2) -- (2,2);
				\draw (0,2) -- (3,0)--(3,2);
				\draw (0,2) -- (0,0);
				\draw (0,2) -- (2,2);
				\draw (0,2) -- (2,0);
				\draw (2,2) -- (0,0);
				\draw (2,2) -- (2,0);
				\draw [fill] (0,0) circle [radius=0.1];
				\node [below] at (0,-0.1) {$3$};
				\draw [fill]  (2,0) circle [radius=0.1];
				\node [below] at (2,-0.1) {$4$};
				\draw [fill]  (0,2) circle [radius=0.1];
				\node [above] at (0,2.1) {$1$};
				\draw [fill]  (2,2) circle [radius=0.1];
				\node [above] at (2,2.1) {$2$};
				\draw [fill]  (3,2) circle [radius=0.1];
				\node [above] at (3,2.1) {$5$};
				\draw [fill]  (3,0) circle [radius=0.1];
				\node [below] at (3,-0.1) {$6$};
			\end{tikzpicture}
		\end{center}
		\caption{The Graph $G(w)$ associated with the word $w$.}
		\label{fig:results}
	\end{figure}
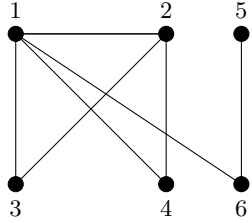
	
\end{example}

For simplicity of arguments, from now on, we mainly 
concentrate on the alphabet with three letters. 
\\
Recall that by a non-trivial component of 
a simple graph, we mean a component which contains 
\emph{at least} one edge.  

\begin{definition}
	For two positive integers $r,s (r > s)$, a graph 
	$G$ is called $(r,s)$-\emph{chordal} graph if 
	for each cycle of $G$ of 
	length at least $r$, there exists at least $s$ chords.  
\end{definition} 

\begin{definition}
	For the given alphabet $\mathcal{A}$,
	The set of all words over the alphabet 
	$\mathcal{A}$ will be denoted by $\mathcal{A}^{\star}$. 
	We also denote the set of non-empty words by 
	$\mathcal{A}^{+}$.
	A word $u$ is called a factor 
	of a word $w$, if there exists 
	words $w_{1}$ and $w_{2}$ such that 
	$w = w_{1} u w_{2}$. 
	A word $u$ is called a prefix (respectively, suffix) 
	of a word $w$, if there exists 
	a words $w_{2}$ 
	(respect., $w_{1}$) such that 
	$w= u w_{2}$ (resp., $w=w_{1}u$).
	A factor $u$ of $w$ which is not equal to $w$ is 
	called a \emph{proper}
	factor. 
\end{definition}

\begin{definition}
	The $k$-th power of a word $w$ is defined 
	\emph{recursively} by 
	$w^{k} = ww^{k-1}$ with the convention that 
	$w^{0} =\lambda$. 
	A word $w \in \mathcal{A}^{+}$ is called 
	\emph{primitive} if the equation $w = u^{n}$ 
	($ u \in \mathcal{A}^{+}$) implies $n=1$. 
	Two words $w$ and $u$ are \emph{conjugate} 
	if there exist two words $w_{1}$
	and $w_{2}$ such that $w = w_{1} w_{2}$
	and $u = w_{2} w_{1}$. 
	It is easy to see that the conjugacy relation is 
	an \emph{equivalence relation}. A \emph{conjugacy class} 
	(or \emph{necklace})
	is a class of this equivalence relation. 
\end{definition}

\begin{definition}
	A word is called a \emph{Lyndon} word 
	if it is primitive and the \emph{smallest} word with respect to 
	the lexicographic order in its conjugacy
	class.
\end{definition}

For now on, we will denote the set of all Lyndon words 
over the alphabet $\mathcal{A}$ of $n$-letters by $Lyn_{n}$

\begin{theorem}
	Any word $w\in \mathcal{A}^{+}$ can be written uniquely as a
	non-increasing product of Lyndon words:
	$$
	w=l_{1}l_{2}\cdots l_{h},~l_{i}\in Lyn_n,
	\hspace{0.5cm}l_{1} 
	\succeq l_{2} \succeq 
	\cdots \succeq l_{h}.
	$$
\end{theorem}

\begin{proposition}
	A word $w\in \mathcal{A}^{+}$ is a Lyndon word 
	if and only if $w\in \mathcal{A}$ or
	$w=rs$ with $r,s\in Lyn_n$ and 
	$r \prec s$. 
	Moreover, if there
	exists a pair $(r, s)$ with $w = rs$ such that $s,w\in
	Lyn_{n}$ and $s$ of maximal length, then 
	$r\in Lyn_{n}$
	and $r \prec rs \prec s$.
\end{proposition}

\subsection{
	Main Results 
}

The following results have similar proofs as in reference 
\cite{BM2016}. Thus, we include them without proofs. 

\begin{lemma}
	A disconnected 
	graph having more than two non-trivial 
	component is not Lyndon-word word representable for any word 
	over ternary alphabet. 
\end{lemma}

\begin{lemma}
	A connected Lyndon-word representable graph 
	with $N$ vertices must have two adjacent vertices whose degree sum is $N$.
\end{lemma}

\begin{proposition}
	Any Lyndon-word representable graph $G=G(w)$ over an ordered 
	ternary alphabet 
	$\mathcal{A}=\{
	a \prec b \prec c 
	\}$
	is a tripartite graph with tripartition 
	$(A,B,C)$, as follows
	$$
	A=\{
	i \in V(G) ~\vert~ w_{i} = a
	\}, \hspace{0.1cm}
	B=\{
	j \in V(G) ~\vert~ w_{j} = b
	\}, \hspace{0.1cm}
	C=\{
	k \in V(G) ~\vert~ w_{k} = c
	\}.
	$$
\end{proposition}

The above proposition immediately implies the 
following results.

\begin{corollary}
	If an \emph{induced subgraph}
	$H$ of a graph $G$ is not word 
	representable, then the original graph 
	$G$ is not word representable. 
\end{corollary}

Finally, we present the main result of this section
which is a characterization of the class of ternary 
Lyndon-word representable graphs. 

\begin{theorem}
	A connected tripartite graph $G$ is 
	ternary Lyndon-word representable if and only if 
	$G$ is $(6,3)$-chordal having two adjacent vertices whose degree sum is the number of vertices of $G$.  
\end{theorem}

\begin{lemma}\label{keylynd1}
	Let $l \in Lyn_{k}$ be a Lyndon word. 
	Then, the graph $G=G(l)$
	that represents the word $l$
	is always a connected graph. 
\end{lemma}

\begin{lemma}
	Let $l_{r}, l_{s} \in Lyn_{k}$ be two Lyndon 
	words such that $l_{r} \preceq l_{s}$. Put
	$l=l_{s} l_{r}$. 
	Then, the graph $G=G(l)$
	that represents the word $l$
	is also a connected graph. 
\end{lemma}

\begin{proposition}
	Let 
	$
	w=l_{1}l_{2}\cdots l_{h},~l_{i}\in Lyn_n,
	\hspace{0.5cm}l_{1} 
	\succ l_{2} \succ 
	\cdots \succ l_{h}
	$
	,
	be the Lyndon factorization of the word $w$ 
	into $h$ \emph{distinct} Lydon words
	over a ternary alphabet. Then, we have 
	\begin{equation}\label{keycomp1}
		h \geq comp(G(w)),
	\end{equation}
	where $com(G)$ denotes the number of connected components
	of $G$. 
\end{proposition}

\subsection{
	Open Problems and Conjectures
}

In this section, we conclude the talk with several 
interesting open questions and conjectures. 
\\
We define the \emph{generalized Stirling cycle number}
$S(N;n,k)$ as the number of distinct Lyndon factorization 
of the word $w$ of length $N$ into $k$ distinct Lyndon words
over the alphabet $\mathcal{A}$ of $n$-letters. 
By convension, we define $S(N,n,0)=1$.

\begin{conjecture}[Coin Arrangments Identity]
	For any $N\geq n>1$, we have
	\begin{equation*}\label{keycoin1}
		\sum_{k=0}^{N}(-1)^{k} S(N;n,k) = 0.
	\end{equation*}
\end{conjecture}

\begin{question}
	Let $G=G(w)$ be a Lyndon-word representable graph with
	$
	w=l_{1}l_{2}\cdots l_{h},~l_{i}\in Lyn_{n},
	\hspace{0.2cm}l_{1} 
	\succeq l_{2} \succeq 
	\cdots \succeq l_{h}
	$
	.
	Does it true that any connected component 
	$G_{con}(w_{con})$ of $G(w)$ is a Lyndon-word
	representable graph with 
	$w_{con}=\prod_{j\in S} l_{j}$, where 
	$S=\{j_1 \leq j_2 \leq \cdots \leq j_l\}$ with $l \leq h$? 
	
\end{question}

Motivated by the definition of the (classic) \emph{Stirling cycle polynomial}, we 
define the \emph{generalized Stirling cycle polynomial}
by 
$
S_{N}(x) =
\sum_{k=0}^{N}S(N;n,k) x^{k}  
$ 
\begin{conjecture}
	For any $N \geq n \geq 1$, the 	
	generalized Stirling cycle polynomial $S_{N}(x)$ 
	has only real roots. In particular, $x=-1$ is always 
	a real root of $S_{N}(x)$.  
\end{conjecture}

We call a word $w$ an $l$-\emph{uniform} word, if 
any letter of alphabet appears in $w$ exactly $l$ times in $w$. Note that $l=1$, corresponds for the case of permutations.   
\begin{question}
	Can we find a recurrence relation for the generalized 
	Stirling cycle numbers $S(N;n,k)$, whenever our 
	words are $l$-uniform? 
\end{question}

\begin{question}
	Can we give a bijective proof of 
	the coin arrangments identity in 
	Conjecture \ref{keycoin1}?
\end{question}

\begin{question}
	For which classes of graphs the equality holds in 
	inequality \ref{keycomp1}? 
\end{question}




\begin{thebibliography}{99}


\bibitem{BM2016} 
Bera, S. and Mahalingam, K. {\it Structural Properties of Word Representable Graphs}, Mathematics in Computer Science, {\bf 10 },  209-222, 2016. 



\end{thebibliography}
\end{document}